\documentclass[12pt]{amsart}

\usepackage{amssymb}
\usepackage{amsmath}
\usepackage{amsthm,cite}

\makeatletter \@addtoreset{equation}{section} \makeatother

\numberwithin{equation}{section}

\newcommand\lam{\lambda}

\newcommand\Rc{\textup{Rc}}

\begin{document}

\bibliographystyle{plain}
\title[]
{Quasi-convergence of the Ricci flow on locally homogeneous closed 4-manifolds}

\author{Songbo Hou}
\address{Department of Applied Mathematics, College of Science, China Agricultural
University,  Beijing, 100083, P.R. China}
\email{housb10@163.com}

\date{}
\def\baselinestretch{1}

\begin{abstract}
We study the quasi-convergence equivalence of some families of metrics on locally homogeneous closed 4-manifolds with trivial isotropy group, and identify the dimension of each equivalence class under certain conditions.

\end{abstract}
\maketitle \baselineskip 18pt
\section{Introduction}

The Ricci flow  is  a  parabolic partial differential equation system:

$$
\frac{\partial g_{ij}}{\partial t}=-2\Rc, \,\,\,g(0)=g_{0}, \eqno (1)
$$
which plays an important role in studying geometries and topologies of manifolds.

Since homogeneous 3-manifolds are  models of Thurston's  geometrization conjecture, it is natural and essential to study  long-time behaviors on closed three-manifolds. In \cite{IJ92}, J. Isenberg and M. Jackson studied  and described characteristic behaviors of Ricci flow in every class of locally homogeneous  geometries.

The results in \cite{Ham86, Ham97} indicated that the Ricci flow could be useful to study geometric and topological properties of 4-manifolds.  In order to
explore further possibility in four dimensions, J. Isenberg, M. Jackson and Peng
Lu \cite{IJL06} studied the Ricci flow on locally homogeneous four-manifolds, and  found that the Ricci flow has the similar  behaviors as those in \cite{IJ92},  in general, if a solution exists for all time, then the flow has a type III singularity in the sense of Hamilton.

 Recall that the normalized Ricci flow
$$\frac{\partial g_{ij}}{\partial t}=-2\Rc+\frac{2r}{n}g \eqno (2)$$ where $n$ is the dimension of the manifold and $r$ is the average of the scalar curvature $R$. The normalization keeps the volume  constant under the flow. In recent years, there are some sdudies on the Ricci flow under more sophisticated procedures on locally homogeneous manifolds. In \cite{CSC09}, Xiaodong Cao and Laurent Saloff-Coste studied the backward Ricci flow, i.e., assuming  the existence interval of the solution to (2) on locally homogeneous 3-manifolds  is $(-T_{b}, T_{f})$, they described the behavior of the Ricci flow as $t$ goes to $-T_{b}$ and obtained the convergence to a sub-Riemannian geometry by a proper re-scaling. Similar results hold for some classes of locally homogeneous 4-manifolds \cite{CHL}. Related references include \cite{CGS09,CNS08,CSC08,GD08, JL07}.

 As the Ricci flow can only converge to Einstein metrics, which means that the right-hand side of equation (2) becomes zero, there are many examples show that the Ricci flow  does not converge although it exists for all positive time. Many of those examples collapse in the sense that the maximum injectivity radius of the solution to equation (2) goes to $0$ as $t\rightarrow \infty$. In \cite{HI93}, Hamilton and Isenberg used the concept of quasi-convergence to describe the behavior of Ricci flow of a family of solv-geometry metrics on twisted torus bundles. The Ricci flow collapse in this family. Further, in 2000, Knopf \cite{DK00} proved that the quasi-convergence equivalence of any metric in this family contains a 1-parameter family of locally homogeneous metrics,  and he formulated the following definition.
\vspace{1ex}

\noindent{\bf Definition 1.1.} {\it If $g$ and $h$ are evolving Riemannian metrics on $M$, $g$ is said to quasi-convergence to $h$ and denote by $g\in [h]$ if for any $\epsilon>0$ there exists a time $t_{\epsilon}$ such that $$\sup_{M\times [t_{\epsilon},\infty)}|g-h|_{h}<\epsilon.$$}

\vspace{1ex}
In \cite{DK01}, Dan Knopf and Kevin McLeod studied the quasi-convergence of all locally homogeneous metrics for which the Ricci flow exists to infinity. On any locally homogeneous 3-maniflod, there exists a Milnor frame \cite{JM76} which can simultaneously diagonalize the initial metric and the Ricci tensor. Then the Ricci flow reduces to an ODE system. Knopf-McLeod  analyzed $[g]$
 in two case: the diagonal cases and the general cases. In the diagonal case, they assumed that $g$ and $h$ are diagonal under the same Milnor frame. Correspondingly, in the general case, $g$ and $h$ are diagonal under  two different frames. The later case needs more analysis.

 The subject of this paper is to study the quasi-convergence of the Ricci flow on locally homogeneous closed 4-manifolds. A class of four dimensional homogeneous geometries can be identified by $(M,G,I)$ where $M$ is a simply connected four manifolds, $G$ is a transitive Lie group acting on $M$ and $I$ is the minimal isotropy group of the action. Four dimensional homogeneous geometries can be divided into two categories. One category with trivial isotropy group is labelled by $A$. Another category with non-trivial isotropy group is labelled by $B$ (refer to \cite{IJL06}). The locally homogeneous 4-manifolds are quite different from 3-manifolds. In the case of four dimensions, even if we can diagonalize some locally homogeneous metrics, the Ricci flow will destroy the diagonalization as time goes on. In order to overcome such difficulty, J. Isenberg, M. Jackson and P. Lu \cite{IJL06} used transition matrices to identify some families of the initial metrics so that the Ricci flow keeps the diagonalization. Then the Ricci flow equation becomes an ODE system. Motivated by the ideas in \cite{DK00, DK01, IJL06}, we consider two cases when deal with the quasi-convergence on locally homogeneous 4-manifolds. The first case is that, assuming we can use the same  matrices as in \cite {IJL06} to diagonalize initial metrics $g_{0}$ and $\bar{g}_{0}$ with which $g(t)$ and $\bar{g}(t)$ evolves by the Ricci flow, compare $g(t)$ and $\bar{g}(t)$. In the second case, under the assumption that we can use the same types of  matrices with different entries to diagonalize initial metrics $g_{0}$ and $\bar{g}_{0}$ such that $g(t)$ and $\bar{g}(t)$ keep diagonalization, compare $g(t)$ and $\bar{g}(t)$ to study $[g]$.

Denote $\alpha=(Y_{1},Y_{2},Y_{3},Y_{4})$  a frame, and $\mathcal{D_{\alpha}}$ the set of all metrics and their Ricci tensors being diagonal with respect to $\alpha$. We always begin our analysis of $[g]$ by studying

$$[g]_{\alpha}\doteqdot [g]\cap \mathcal{D_{\alpha}}.$$

Denote  $\beta=(Y_{1}^{'},Y_{2}^{'},Y_{3}^{'},Y_{4}^{'})$  another frame. Assume the Ricci flow solution $g(t)$ is diagonal under $\alpha$, and $\bar{g}(t)$ is diagonal under $\beta$. If the transformation from $\beta$ to $\alpha$ is

$$\left[  \begin{array}{c}
  Y_{1}\\
   Y_{2} \\
   Y_{3}\\
  Y_{4}\\
   \end{array}\right]=A\left[  \begin{array}{c}
  Y^{'}_{1}\\
   Y^{'}_{2} \\
   Y^{'}_{3}\\
  Y^{'}_{4}\\
   \end{array}\right]$$ where $A$ is a transition matrix,    the following relation  between evolving component matrix $\bar{g}_{\alpha}(t)$ of $\bar{g}(t)$  under $\alpha$ and $\bar{g}_{\beta}(t)$ under $\beta$ holds:
 $$\bar{g}_{\alpha}(t)=A\bar{g}_{\beta}(t)A^{T}.$$ Then computation on   $|g-\bar{g}|^{2}_{g}$ yields the quasi-convergence result.

  In this paper, we use the same notations as those in \cite{IJL06} to denote the classification of homogeneous 4-mianofolds, and also broaden the discussion to Ricci flow on 4-dimensional unimodular Lie groups and focus on category A including class $\bf{A1}-\bf{A10}$.  Since the existence interval of the Ricci flow in class $\bf{A10}$ is finite, we omit this class. Denoting $X_{1},\,\,X_{2},\,\,X_{3},\,\,X_{4}$ a basis such that the  Lie brackets take the characteristic forms ( refer to \cite {IJL06,MM92}), and $$g_{0}=\lambda_{1}\omega^{1}\otimes\omega^{1}+\lambda_{2}\omega^{2}\otimes\omega^{2}+\lambda_{3}\omega^{3}\otimes\omega^{3}+\lambda_{4}\omega^{4}\otimes\omega^{4},$$
  where $\omega^{i}$ is the frame of 1-forms dual to $X_{i}$, we write
  $$g(t)=A(t)\omega^{1}\otimes\omega^{1}+B(t)\omega^{2}\otimes\omega^{2}+C(t)\omega^{3}\otimes\omega^{3}+D(t)\omega^{4}\otimes\omega^{4}.$$
  In every section, we   recall  the Ricci flow equation, the solution or the long-time behavior, and the conserved quantities in \cite{IJL06}, then analyze every case. We will get the conditions  some entries in transition matrices satisfy and the dimension of each quasi-convergence equivalence class.

\section {\bf {A1.} $U1[(1,1,1)]$}

The Ricci flow is trivial and $[g]$ is $0$-dimensional.

\section {\bf {A2.}  $U1[1,1,1].$}
\noindent{\bf 3.1.  The class $[g]_{\alpha}$}

\vspace{1ex}
\noindent{\bf (A2iv)} This corresponds to case in which $k\neq0,\,1,\,-\frac{1}{2}$. From the conclusions in \cite{IJL06}, assume that we can use $Y_{i}=\Lambda_{i}^{k}X_{k}$ with

 $$\Lambda=\left[  \begin{array}{cccc}
  1 & 0 & 0 & 0 \\
   0 & 1 & 0 & 0 \\
   0 & 0 & 1 &0 \\
  a_{4} & a_{5} & a_{6} & 1 \\
   \end{array}\right]$$ to diagonalize the initial metric $g_{0}$.
  The Ricci flow is
\begin{align*}
&\frac{dA}{dt}=0,\,\,&\frac{dB}{dt}=&0,\\
&\frac{dC}{dt}=0,\,\, &\frac{dD}{dt}=&4(k^{2}+k+1).
\end{align*}
It is easy to obtain the solution
\begin{align*}
&A(t)=\lambda_{1},\,\,&B(t)=&\lambda_{2},\\
&C(t)=\lambda_{3},\,\, &D(t)=&\lambda_{4}+4(k^{2}+k+1)t.
\end{align*}
Denote $\alpha=(Y_{1},Y_{2},Y_{3},Y_{4})$.
If $$\bar{g}=\bar{A}\omega^{1}\otimes\omega^{1}+\bar{B}\omega^{2}\otimes\omega^{2}+\bar{C}\omega^{3}\otimes\omega^{3}+\bar{D}\omega^{4}\otimes\omega^{4}$$ is another metric  in $[g]_{\alpha}$ with initial data $(\bar{\lambda}_{1},\,\bar{\lambda}_{2},\,\bar{\lambda}_{3},\,\bar{\lambda}_{4})$, then
$$|g-\bar{g}|^{2}_{g}=\left(\frac{A-\bar{A}}{A}\right)^{2}+\left(\frac{B-\bar{B}}{B}\right)^{2}+\left(\frac{C-\bar{C}}{C}\right)^{2}+\left(\frac{D-\bar{D}}{D}\right)^{2},$$
where
$$\frac{A-\bar{A}}{A}=1-\frac{\bar{\lambda}_{1}}{\lambda_{1}},$$
$$\frac{B-\bar{B}}{B}=1-\frac{\bar{\lambda}_{2}}{\lambda_{2}},$$
$$\frac{C-\bar{C}}{C}=1-\frac{\bar{\lambda}_{3}}{\lambda_{3}},$$
$$\frac{D-\bar{D}}{D}=1-\frac{\bar{\lambda}_{4}+4(k^{2}+k+1)t}{\lambda_{4}+4(k^{2}+k+1)t}.$$

\vspace{1ex}
\noindent {\bf Lemma 3.1.} {\it The class $[g]_{\alpha}$ for an {\bf A2iv}-geometry metric is exactly a 1-parameter family.}

\vspace{1ex}
\noindent {\bf Proof.} Note that $\bar{g}\in [g]_{\alpha}$ iff $\frac{\bar{A}}{A}$, $\frac{\bar{B}}{B}$, $\frac{\bar{C}}{C}$, and $\frac{\bar{D}}{D}$ all converge to 1.
It follows that  $\bar{\lambda}_{1}=\lambda_{1}$, $\bar{\lambda}_{2}=\lambda_{2}$, $\bar{\lambda}_{3}=\lambda_{3}$ and $\bar{\lambda}_{4}$ can be chosen arbitrarily. The proof is got.

\vspace{1ex}
\noindent{\bf 3.2.  $g$ and $\bar{g}$ are diagonal under two different frames $\alpha$ and $\beta$}
\vspace{1ex}

Assume that we can use  $Y_{i}=\Lambda_{i}^{k}X_{k}$ with

 $$\Lambda=\left[  \begin{array}{cccc}
  1 & 0 & 0 & 0 \\
   0 & 1 & 0 & 0 \\
   0 & 0 & 1 & 0\\
  a_{4} & a_{5} & a_{6} & 1 \\
   \end{array}\right]$$ and $Y^{'}_{i}=\Lambda_{i}^{'k}X_{k}$ with

 $$\Lambda^{'}=\left[  \begin{array}{cccc}
  1 & 0 & 0 & 0 \\
   0 & 1 & 0 & 0 \\
   0 & 0 & 1 & 0\\
  a_{4}^{'} & a_{5}^{'} & a_{6}^{'} & 1 \\
   \end{array}\right]$$to diagonalize $g_{0}$,  $\bar{g}_{0}$ respectively.

   Denote  $\beta=(Y_{1}^{'},Y_{2}^{'},Y_{3}^{'},Y_{4}^{'})$ . The transformation for $\mathcal{\beta}$ to $\alpha$ is
   $$\left[  \begin{array}{c}
  Y_{1}\\
   Y_{2} \\
   Y_{3}\\
  Y_{4}\\
   \end{array}\right]=\Lambda\Lambda^{'-1} \left[  \begin{array}{c}
  Y^{'}_{1}\\
   Y^{'}_{2} \\
   Y^{'}_{3}\\
  Y^{'}_{4}\\
   \end{array}\right].$$
   Letting $$a=a_{4}-a_{4}^{'},\,\,b=a_{5}-a_{5}^{'},\,\,c=a_{6}-a_{6}^{'},$$ and
   $$A=\left[  \begin{array}{cccc}
  1 & 0 & 0 & 0 \\
   0 & 1 & 0 & 0 \\
   0 & 0 & 1 & 0\\
  a & b & c & 1 \\
   \end{array}\right],$$  we have
   $$\left[  \begin{array}{c}
  Y_{1}\\
   Y_{2} \\
   Y_{3}\\
  Y_{4}\\
   \end{array}\right]=A\left[  \begin{array}{c}
  Y^{'}_{1}\\
   Y^{'}_{2} \\
   Y^{'}_{3}\\
  Y^{'}_{4}\\
   \end{array}\right].$$
 The relation between $\bar{g}_{\alpha}(t)$ and $\bar{g}_{\beta}(t)$ is
 $$\bar{g}_{\alpha}(t)=A\bar{g}_{\beta}(t)A^{T}.$$
 Assuming  $$\bar{g}_{\beta}(t)=\left[  \begin{array}{cccc}
  \bar{A} & 0 & 0 & 0 \\
   0 & \bar{B} & 0 & 0 \\
   0 & 0 & \bar{C} & 0\\
  0 & 0 & 0 & \bar{D} \\
   \end{array}\right],$$
then we get $$\bar{g}_{\alpha}(t)=\left[  \begin{array}{cccc}
  \bar{A} & 0 & 0 & \bar{A}a \\
   0 & \bar{B} & 0 & \bar{B}b \\
   0 & 0 & \bar{C} & \bar{C}c\\
  \bar{A}a & \bar{B}b & \bar{C}c & a^{2}\bar{A}+b^{2}\bar{B}+c^{2}\bar{C}+\bar{D} \\
   \end{array}\right].$$
\vspace{1ex}

\noindent {\bf Theorem 3.1. } {\it The quasi-convergence class $[g]$ is exactly a 1-parameter family.}

\vspace{1ex}
 \noindent{\bf Proof.} In fact, it is easy to see that
 \begin{align*}
&\bar{A}(t)=\bar{\lambda}_{1},\,\,&\bar{B}(t)=&\bar{\lambda}_{2},\\
&\bar{C}(t)=\bar{\lambda}_{3},\,\, &\bar{D}(t)=&\bar{\lambda}_{4}+4(k^{2}+k+1)t.
\end{align*}
Thus $\bar{g}\in [g]$ if and only if all terms convergence to 0 in the sum
 \begin{align*}|
\bar{g}-g|_{g}^{2}=&\left(\frac{A-\bar{A}}{A}\right)^{2}+\left(\frac{B-\bar{B}}{B}\right)^{2}+\left(\frac{C-\bar{C}}{C}\right)^{2}+2\frac{(\bar{A}a)^{2}}{AD}\\
  +&\frac{2(\bar{B}b)^{2}}{BD}+\frac{2(\bar{C}c)^{2}}{CD}+\left(\frac{a^{2}\bar{A}+b^{2}\bar{B}+c^{2}\bar{C}+\bar{D}-D}{D}\right)^{2}.\\
 \end{align*}
  This implies that $\bar{\lambda}_{1}$, $\bar{\lambda}_{2}$ and $\bar{\lambda}_{3}$ are determined, and $\bar{\lambda}_{4}$ is arbitrary.

 The analysis in the cases {\bf (A2i)}, {\bf (A2ii)}, {\bf (A2iii)} is similar.

\section {\bf {A3.}  $U1[Z,\bar{Z},1].$}

\noindent{\bf 4.1. The class $[g]_{\alpha}$ }
\vspace{1ex}

 From the conclusions in \cite{IJL06},
assume that we can use $Y_{i}=\Lambda_{i}^{k}X_{k}$ with

 $$\Lambda=\left[  \begin{array}{cccc}
  1 & 0 & 0 & 0 \\
   0 & 1 & 0 & 0 \\
   0 & 0 & 1 & 0\\
  a_{4} & a_{5} & a_{6} & 1 \\
   \end{array}\right]$$ to diagonalize the initial metric $g_{0}$.

   The Ricci flow is
\begin{align*}
&\frac{dA}{dt} =-\frac{A^{2}-B^{2}}{BD}, &\frac{dB}{dt} =&-\frac{B^{2}-A^{2}}{AD}, \\
&\frac{dC}{dt} =0, &\frac{dD}{dt} =&\frac{(A-B)^{2}+12k^{2}AB}{AB}.\end{align*}

   From the behavior of the solution to Ricci flow, we get the following theorem.

\vspace{1ex}
\noindent {\bf Lemma 4.1.} {\it The class $[g]_{\alpha}$  is exactly a 2-parameter family.}

\vspace{1ex}
\noindent{\bf Proof.}   If $\lambda_{1}=\lambda_{2}$,
    then $$A(t)=\lambda_{1},\,\,B(t)=\lambda_{2},\,\,\,\,C(t)=\lambda_{3},\,\,D(t)=\lambda_{4}+12k^{2}t.$$
    If $\lambda_{1}\neq\lambda_{2}$,the long-time behavior of the solution $g(t)$ is (refer to \cite{IJL06})
$$A(t)\rightarrow \sqrt{\lambda_{1}\lambda_{2}},\,\,B(t)\rightarrow \sqrt{\lambda_{1}\lambda_{2}},\,\,C(t)=\lambda_{3},\,\, D(t)\rightarrow \infty.$$
Let $\bar{g}\in[g]_{\alpha}$ correspond to the initial data $(\bar{\lambda}_{1},\bar{\lambda}_{2},\bar{\lambda}_{3},\bar{\lambda}_{4})$ . Then  $\frac{\bar{A}}{A}\rightarrow 1$ and $\frac{\bar{B}}{B}\rightarrow 1$ if and only if $\bar{\lambda}_{1}\bar{\lambda}_{2}=\lambda_{1}\lambda_{2}$. It is easy to see that
$$\lim_{t\rightarrow \infty}\frac{\bar{D}}{D}=\lim_{t\rightarrow \infty}\frac{(\bar{A}-\bar{B})^{2}+12k^{2}\bar{A}\bar{B}}{(A-B)^{2}+12k^{2}AB}=1.$$ If $\bar{\lambda}_{1}$ and $\bar{\lambda}_{4}$ are prescribed, $\bar{\lambda}_{2}$  and $\bar{\lambda}_{3}$ are determined.

\vspace{1ex}

\noindent{\bf 4.2.  $g$ and $\bar{g}$ are diagonal under two different frames $\alpha$ and $\beta$}

\vspace{1ex}

The analysis in this case is similar to that in Section 3.2.


\section {\bf {A4.}  $U1[2,1].$}

\noindent{\bf 5.1. The class $[g]_{\alpha}$ }
\vspace{1ex}

 From the conclusions in \cite{IJL06}, assume that we can use $Y_{i}=\Lambda_{i}^{k}X_{k}$ with

 $$\Lambda=\left[  \begin{array}{cccc}
  1 & a_{2} & a_{3} & 0 \\
   0 & 1 & 0 & 0 \\
   0 & a_{1} & 1 & 0\\
  a_{4} & a_{5} & a_{6} & 1 \\
   \end{array}\right]$$ to diagonalize the initial metric $g_{0}$.
 The solution to the Ricci flow is
\begin{align*}
&A =\lam_{1}\left(1+\frac{3\lam_{2}}{\lam_{1}\lam_{4}}t\right)^{1/3}, &B =&\lam_{2}\left(1+\frac{3\lam_{2}}{\lam_{1}\lam_{4}}t\right)^{-1/3}, \\
&C =\lam_{3}, &D=&\lam_{4}\left(1+\frac{3\lam_{2}}{\lam_{1}\lam_{4}}t\right)^{1/3}.
\end{align*}

\vspace{1ex}

\noindent {\bf Lemma 5.1.} {\it The class $[g]_{\alpha}$ for an  {\bf A4}-geometry metric is exactly a 1-parameter family.}
\vspace{1ex}

\noindent
{\bf Proof.}
 Assume that $$\bar{g}(t)=\bar{A}\omega^{1}\otimes\omega^{1}+\bar{B}\omega^{2}\otimes\omega^{2}+\bar{C}\omega^{3}\otimes\omega^{3}+\bar{D}\omega^{4}\otimes\omega^{4}.$$

 Then $\bar{g}\in [g]_{\alpha}$ if and only if $\frac{\bar{A}}{A}\rightarrow 1$,$\frac{\bar{B}}{B}\rightarrow 1$,$\frac{\bar{C}}{C}\rightarrow 1$,$\frac{\bar{D}}{D}\rightarrow 1$. Thus we derive $\bar{\lambda}_{3}=\lambda_{3}$. Let $\frac{\bar{\lambda}_{1}}{\lambda_{1}}=k$. We obtain
 $$\frac{\bar{\lambda}_{2}}{\lambda_{2}}=\frac{1}{k}\,\,\,\,\,\,\frac{\bar{\lambda}_{4}}{\lambda_{4}}=k.$$ We choose $\bar{\lambda}_{1}$ arbitrarily. Then $\bar{\lambda}_{2}$ and $\bar{\lambda}_{4}$ will be determined. This completes the proof.
\vspace{1ex}

 \noindent {\bf 5.2.  $g$ and $\bar{g}$ are diagonal under two different frames $\alpha$ and $\beta$}

\vspace{1ex}

Assume that we can use  $Y_{i}=\Lambda_{i}^{k}X_{k}$ with

 $$\Lambda=\left[  \begin{array}{cccc}
  1 & a_{2} & a_{3} & 0 \\
   0 & 1 & 0 & 0 \\
   0 & a_{1} & 1 & 0\\
  a_{4} & a_{5} & a_{6} & 1 \\
   \end{array}\right]$$ and $Y^{'}_{i}=\Lambda_{i}^{'k}X_{k}$ with

 $$\Lambda^{'}=\left[  \begin{array}{cccc}
  1 & a_{2}^{'} & a_{3}^{'} & 0 \\
   0 & 1 & 0 & 0 \\
   0 & a_{1}^{'} & 1 & 0\\
  a_{4}^{'} & a_{5}^{'} & a_{6}^{'} & 1 \\
   \end{array}\right]$$to diagonalize $g_{0}$,  $\bar{g}_{0}$ respectively.

   The transformation for $\mathcal{\beta}$ to $\alpha$ is
   $$\left[  \begin{array}{c}
  Y_{1}\\
   Y_{2} \\
   Y_{3}\\
  Y_{4}\\
   \end{array}\right]=\Lambda\Lambda^{'-1} \left[  \begin{array}{c}
  Y^{'}_{1}\\
   Y^{'}_{2} \\
   Y^{'}_{3}\\
  Y^{'}_{4}\\
   \end{array}\right].$$
   Denoting $$a=a_{1}-a_{1}^{'},\,\,b=a_{2}-a_{2}^{'}+a_{1}^{'}a_{3}^{'}-a_{1}^{'}a_{3},\,\,c=a_{3}-a_{3}^{'},\,\,d=a_{4}-a_{4}^{'},$$
   $$e=a_{5}-a_{5}^{'}+a_{2}^{'}a_{4}^{'}-a_{2}^{'}a_{4}+a_{1}^{'}a_{6}^{'}-a_{1}^{'}a_{6}+a_{1}^{'}a_{3}^{'}a_{4}-a_{1}^{'}a_{3}^{'}a_{4}^{'},$$
   $$f=a_{6}-a_{6}^{'}+a_{3}^{'}a_{4}^{'}-a_{3}^{'}a_{4}$$
    and
   $$A=\left[  \begin{array}{cccc}
  1 & b & c & 0 \\
   0 & 1 & 0 & 0 \\
   0 & a & 1 & 0\\
  d & e & f & 1 \\
   \end{array}\right],$$  we have
   $$\left[  \begin{array}{c}
  Y_{1}\\
   Y_{2} \\
   Y_{3}\\
  Y_{4}\\
   \end{array}\right]=A\left[  \begin{array}{c}
  Y^{'}_{1}\\
   Y^{'}_{2} \\
   Y^{'}_{3}\\
  Y^{'}_{4}\\
   \end{array}\right].$$
  The relation between $\bar{g}_{\alpha}(t)$ and $\bar{g}_{\beta}(t)$ is
 $$\bar{g}_{\alpha}(t)=A\bar{g}_{\beta}(t)A^{T}.$$
 Assuming $$\bar{g}_{\beta}(t)=\left[  \begin{array}{cccc}
  \bar{A} & 0 & 0 & 0 \\
   0 & \bar{B} & 0 & 0 \\
   0 & 0 & \bar{C} & 0\\
  0 & 0 & 0 & \bar{D} \\
   \end{array}\right],$$
 then we get $$\bar{g}_{\alpha}(t)= \left[ \begin {array}{cccc} \bar{A}+{b}^{2}\bar{B}+{c}^{2}\bar{C}&b\bar{B}&b\bar{B}a+c\bar{C}&d\bar{A}+b\bar{B}e+c\bar{C}f\\ \noalign{\medskip}b\bar{B}&\bar{B}&a\bar{B}&e\bar{B}\\ \noalign{\medskip}b\bar{B}a+c\bar{C}&
a\bar{B}&{a}^{2}\bar{B}+\bar{C}&a\bar{B}e+f\bar{C}\\ \noalign{\medskip}d\bar{A}+b\bar{B}e+c\bar{C}f&e\bar{B}&a\bar{B}e+f\bar{C}&{d
}^{2}\bar{A}+{e}^{2}\bar{B}+{f}^{2}\bar{C}+\bar{D}\end {array} \right]
.$$
 Thus $\bar{g}\in [g]$ if and only if all terms convergence to 0 in the sum

\begin{align*}|
\bar{g}-g|_{g}^{2}=&\left(\frac{\bar{A}+b^{2}\bar{B}+c^{2}\bar{C}-A}{A}\right)^{2}+\left(\frac{B-\bar{B}}{B}\right)^{2}+\left(\frac{C-\bar{C}-a^{2}\bar{B}}{C}\right)^{2}\\
+&\left(\frac{d^{2}\bar{A}+e^{2}\bar{B}+f^{2}\bar{C}+\bar{D}-D}{D}\right)^{2}+2\frac{(b\bar{B})^{2}}{AB}+2\frac{(b\bar{B}a+c\bar{C})^{2}}{AC}+2\frac{(d\bar{A}+b\bar{B}e+c\bar{C}f)^{2}}{AD}\\
+2&\frac{(a\bar{B})^{2}}{BC}+2\frac{(e\bar{B})^{2}}{BD}+2\frac{(a\bar{B}e+f\bar{C})^{2}}{CD}.\\
  \end{align*}
  It follows that  $d=0$ which  implies that $a_{4}^{'}=a_{4}$, and $\frac{\bar{A}}{A}\rightarrow 1$,$\frac{\bar{B}}{B}\rightarrow 1$,$\frac{\bar{C}}{C}\rightarrow 1$,$\frac{\bar{D}}{D}\rightarrow 1$.  Under the condition $a_{4}=a_{4}^{'}$, we have the following Theorem.

\vspace{1ex}

\noindent {\bf Theorem 5.1.} {\it The class $[g]$ is  exactly  a 1-parameter family.}

\vspace{1ex}
The proof is similar to Lemma 5.1.


\section {\bf {A5.}  $U1[2,1].$}
\vspace{1ex}

\noindent{\bf 6.1. The class $[g]_{\alpha}$ }
\vspace{1ex}

Recall  the results in \cite{IJL06}. The Ricci flow is
\begin{align*}
&\frac{dA}{dt} =\frac{B}{D}=\frac{B}{AD}A, &\frac{dB}{dt} =&-\frac{B^{2}}{AD}=-\frac{B}{AD}B, \\
&\frac{dC}{dt} =0, &\frac{dD}{dt} =&3+\frac{B}{A}.
\end{align*}
The  long time behavior of the solution $g(t)$ is
$$A(t)\rightarrow \infty,\,\,B(t)\rightarrow 0^{+},\,\, C(t)=\lambda_{3},\,\,D(t)\rightarrow \infty.$$

We also have $AB=\lambda_{1}\lambda_{2}$.

\vspace{1ex}

\noindent {\bf Lemma 6.1.} {\it The class $[g]_{\alpha}$ for an  {\bf A5}-geometry metric is exactly a 2-parameter family.}
\vspace{1ex}

\noindent {\bf Proof.} Assume that $$\bar{g}(t)=\bar{A}\omega^{1}\otimes\omega^{1}+\bar{B}\omega^{2}\otimes\omega^{2}+\bar{C}\omega^{3}\otimes\omega^{3}+\bar{D}\omega^{4}\otimes\omega^{4}$$ is another Ricci flow solution with the initial data $(\bar{\lambda}_{1},\bar{\lambda}_{2},\bar{\lambda}_{3},\bar{\lambda}_{1})$.
 Then $\bar{g}\in [g]_{\alpha}$ if and only if all terms converge to $0$ in the sum
$$|g-\bar{g}|^{2}_{g}=\left(\frac{A-\bar{A}}{A}\right)^{2}+\left(\frac{B-\bar{B}}{B}\right)^{2}+\left(\frac{C-\bar{C}}{C}\right)^{2}+\left(\frac{D-\bar{D}}{D}\right)^{2}.$$
 From the equation, it follows that
$$\lim_{t\rightarrow\infty}\frac{\bar{D}(t)}{D(t)}=\lim_{t\rightarrow\infty}\frac{\frac{d\bar{D}(t)}{dt}}{\frac{dD(t)}{dt}}
=\lim_{t\rightarrow\infty}\frac{3+\frac{\bar{B}}{\bar{A}}}{3+\frac{B}{A}}=1.$$  Then $\bar{A}/A\rightarrow 1$ and $\bar{B}/B\rightarrow 1$ if and only if $\bar{\lambda}_{1}\bar{\lambda}_{2}=\lambda_{1}\lambda_{2}$ since $\frac{dA^{2}}{dt}=\frac{2\lambda_{1}\lambda_{2}}{D}$.
From $\bar{\lambda}_{2}=\frac{\lambda_{1}\lambda_{2}}{\bar{\lambda}_{1}}$, once $\bar{\lambda}_{1}$ is chosen, $\bar{\lambda}_{2}$ is fixed, and $\bar{\lambda}_{4}$ can be arbitrary.

\vspace{1ex}
 \noindent{\bf 6.2.  $g$ and $\bar{g}$ are diagonal under two different frames $\alpha$ and $\beta$.}
\vspace{1ex}

Assume  we can use  $Y_{i}=\Lambda_{i}^{k}X_{k}$ with

 $$\Lambda=\left[  \begin{array}{cccc}
  1 & a_{2} & 0 & 0 \\
   0 & 1 & 0 & 0 \\
   0 & 0 & 1 & 0\\
  a_{4} & a_{5} & a_{6} & 1 \\
   \end{array}\right]$$ and $Y^{'}_{i}=\Lambda_{i}^{'k}X_{k}$ with

 $$\Lambda^{'}=\left[  \begin{array}{cccc}
  1 & a_{2}^{'} & 0 & 0 \\
   0 & 1 & 0 & 0 \\
   0 & 0 & 1 & 0\\
  a_{4}^{'} & a_{5}^{'} & a_{6}^{'} & 1 \\
   \end{array}\right]$$to diagonalize $g_{0}$,  $\bar{g}_{0}$ respectively.

   The transformation for $\mathcal{\beta}$ to $\alpha$ is
   $$\left[  \begin{array}{c}
  Y_{1}\\
   Y_{2} \\
   Y_{3}\\
  Y_{4}\\
   \end{array}\right]=\Lambda\Lambda^{'-1} \left[  \begin{array}{c}
  Y^{'}_{1}\\
   Y^{'}_{2} \\
   Y^{'}_{3}\\
  Y^{'}_{4}\\
   \end{array}\right].$$
   Denoting $$a=a_{2}-a_{2}^{'},\,\,b=a_{4}-a_{4}^{'},\,\,c=a_{5}-a_{5}^{'}+a_{2}^{'}a_{4}^{'}-a_{2}^{'}a_{4},\,\,d=a_{6}-a_{6}^{'}$$ and
   $$A=\left[  \begin{array}{cccc}
  1 & a & 0 & 0 \\
   0 & 1 & 0 & 0 \\
   0 & 0 & 1 & 0\\
  b & c & d & 1 \\
   \end{array}\right],$$  we have
   $$\left[  \begin{array}{c}
  Y_{1}\\
   Y_{2} \\
   Y_{3}\\
  Y_{4}\\
   \end{array}\right]=A\left[  \begin{array}{c}
  Y^{'}_{1}\\
   Y^{'}_{2} \\
   Y^{'}_{3}\\
  Y^{'}_{4}\\
   \end{array}\right].$$
As previous,
 we get $$\bar{g}_{\alpha}(t)=\left[  \begin{array}{cccc}
  \bar{A}+a^{2}\bar{B} & a\bar{B} & 0 & b\bar{A}+ac\bar{B} \\
   a\bar{B} & \bar{B} & 0 & c\bar{B} \\
   0 & 0 & \bar{C} & d\bar{C}\\
  b\bar{A}+ac\bar{B} & c\bar{B} & d\bar{C} & b^{2}\bar{A}+c^{2}\bar{B}+d^{2}\bar{C}+\bar{D} \\
   \end{array}\right].$$
 Thus $\bar{g}\in [g]$ if and only if all terms convergence to 0 in the sum
\begin{align*}|
\bar{g}-g|_{g}^{2}=&\left(\frac{A-\bar{A}-a^{2}\bar{B}}{A}\right)^{2}+\left(\frac{B-\bar{B}}{B}\right)^{2}+\left(\frac{C-\bar{C}}{C}\right)^{2}+\left(\frac{b^{2}\bar{A}+c^{2}\bar{B}+d^{2}\bar{C}+\bar{D}-D}{D}\right)^{2}\\
  +&\frac{2(a\bar{B})^{2}}{AB}+\frac{2(b\bar{A}+ac\bar{B})^{2}}{AD}+\frac{2(c\bar{B})^{2}}{BD}+\frac{2(d\bar{C})^{2}}{CD}.\\
  \end{align*}It follows that $\bar{\lambda}_{3}=\lambda_{3}$, $\lim_{t\rightarrow 0}\frac{\bar{A}}{A}=1$ and $\lim_{t\rightarrow 0}\frac{\bar{B}}{B}=1$ from $(\frac{C-\bar{C}}{C})^{2}$, $(\frac{A-\bar{A}-a^{2}\bar{B}}{A})^{2}$ and $(\frac{B-\bar{B}}{B})^{2}$  goes to $0$. Other terms go to $0$ obviously after simply computation.
We also have the following Theorem.

\vspace{1ex}
\noindent {\bf Theorem 6.1.} {\it The class $[g]$  is exactly a 2-parameter family.}

\section {\bf {A6.}  $U1[3].$}

\vspace{1ex}

\noindent{\bf 7.1. The class $[g]_{\alpha}$ }
\vspace{1ex}

Recall the results in \cite{IJL06}.
The  Ricci flow is
\begin{align*}
&\frac{dA}{dt} =\frac{B}{D}, &\frac{dB}{dt} =&\frac{AC-B^{2}}{AD}, \\
&\frac{dC}{dt} =-\frac{C^{2}}{BD}, &\frac{dD}{dt} =&\frac{B}{A}+\frac{C}{B}.
\end{align*}
The solution to the Ricci flow is
\begin{align*}
&A(t)=\lambda_{1}(3E_{0}t+1)^{1/3}, &B(t)=\lambda_{2}(3E_{0}t+1)^{-1/3}(3F_{0}t+1)^{1/3}, \\
&C(t)=\lambda_{3}(3F_{0}t+1)^{-1/3}, &D(t)=\lambda_{4}(3E_{0}t+1)^{1/3}(3F_{0}t+1)^{1/3},
\end{align*}
where $E_{0}=\frac{\lambda_{2}}{\lambda_{1}\lambda_{4}}$ and $F_{0}=\frac{\lambda_{3}}{\lambda_{2}\lambda_{4}}$.
\vspace{1ex}

\noindent {\bf Lemma 7.1.} {\it The class $[g]_{\alpha}$ is exactly a 2-parameter family.}
\vspace{1ex}

 \noindent{\bf Proof.} Assume that $$\bar{g}(t)=\bar{A}\omega^{1}\otimes\omega^{1}+\bar{B}\omega^{2}\otimes\omega^{2}+\bar{C}\omega^{3}\otimes\omega^{3}+\bar{D}\omega^{4}\otimes\omega^{4}$$ is another Ricci flow solution with the initial data $(\bar{\lambda}_{1},\bar{\lambda}_{2},\bar{\lambda}_{3},\bar{\lambda}_{4})$.
 Then $\bar{g}\in [g]_{\alpha}$ if and only if all terms converge to $0$ in the sum
$$|g-\bar{g}|^{2}_{g}=\left(\frac{A-\bar{A}}{A}\right)^{2}+\left(\frac{B-\bar{B}}{B}\right)^{2}+\left(\frac{C-\bar{C}}{C}\right)^{2}+\left(\frac{D-\bar{D}}{D}\right)^{2}.$$
Thus $\bar{A}/A\rightarrow 1$,  $\bar{B}/B\rightarrow 1$, $\bar{C}/C\rightarrow 1$ and $\bar{D}/D\rightarrow 1$ if and only if   the following equalities hold
\begin{align*}
&\frac{\bar{\lambda}_{1}}{\lambda_{1}}\left(\frac{\bar{E}_{0}}{E_{0}}\right)^{1/3}=1,&\frac{\bar{\lambda}_{2}}{\lambda_{2}}\left(\frac{E_{0}}{\bar{E}_{0}}\right)^{1/3}\left(\frac{\bar{F}_{0}}{F_{0}}\right)^{1/3}=1,\\
&\frac{\bar{\lambda}_{3}}{\lambda_{3}}\left(\frac{F_{0}}{\bar{F}_{0}}\right)^{1/3}=1, &\frac{\bar{\lambda}_{4}}{\lambda_{4}}\left(\frac{\bar{E}_{0}}{E_{0}}\right)^{1/3}\left(\frac{\bar{F}_{0}}{F_{0}}\right)^{1/3}=1.
\end{align*}
It is easy to see that if any two terms of $\bar{\lambda}_{1},\,\bar{\lambda}_{2},\,\bar{\lambda}_{3},\,\bar{\lambda}_{4}$ are chosen, then the left two will be determined.
\vspace{1ex}

\noindent{\bf 7.2.  $g$ and $\bar{g}$ are diagonal under two different frames $\alpha$ and $\beta$}
\vspace{1ex}

Assume  we can use  $Y_{i}=\Lambda_{i}^{k}X_{k}$ with

 $$\Lambda=\left[  \begin{array}{cccc}
  1 & a_{1} & a_{3} & 0 \\
   0 & 1 & a_{1} & 0 \\
   0 & 0 & 1 & 0\\
  a_{4} & a_{5} & a_{6} & 1 \\
   \end{array}\right]$$ and $Y^{'}_{i}=\Lambda_{i}^{'k}X_{k}$ with

 $$\Lambda^{'}=\left[  \begin{array}{cccc}
  1 & a_{1}^{'} & a_{3}^{'} & 0 \\
   0 & 1 & a_{1}^{'} & 0 \\
   0 & 0 & 1 & 0\\
  a_{4}^{'} & a_{5}^{'} & a_{6}^{'} & 1 \\
   \end{array}\right]$$to diagonalize $g_{0}$,  $\bar{g}_{0}$ respectively.

   The transformation for $\mathcal{\beta}$ to $\alpha$ is
   $$\left[  \begin{array}{c}
  Y_{1}\\
   Y_{2} \\
   Y_{3}\\
  Y_{4}\\
   \end{array}\right]=\Lambda\Lambda^{'-1} \left[  \begin{array}{c}
  Y^{'}_{1}\\
   Y^{'}_{2} \\
   Y^{'}_{3}\\
  Y^{'}_{4}\\
   \end{array}\right].$$
   Denoting $$a=a_{1}-a_{1}^{'},\,\,b=a_{4}-a_{4}^{'},\,\,c=a_{5}-a_{5}^{'}+a_{1}^{'}a_{4}^{'}-a_{1}^{'}a_{4},\,\,
   $$$$d=a_{6}-a_{6}^{'}+a_{1}^{'}a_{5}^{'}-a_{1}^{'}a_{5}+a_{1}^{'2}a_{4}-a_{1}^{'2}a_{4}^{'}+a_{3}^{'}a_{4}^{'}-a_{3}^{'}a_{4},\,\,e=a_{3}-a_{3}^{'}+a_{1}^{'2}-a_{1}a_{1}^{'}.
   $$ and
   $$A=\left[  \begin{array}{cccc}
  1 & a & e & 0 \\
   0 & 1 & a & 0 \\
   0 & 0 & 1 & 0\\
  b & c & d & 1 \\
   \end{array}\right],$$  we have
   $$\left[  \begin{array}{c}
  Y_{1}\\
   Y_{2} \\
   Y_{3}\\
  Y_{4}\\
   \end{array}\right]=A\left[  \begin{array}{c}
  Y^{'}_{1}\\
   Y^{'}_{2} \\
   Y^{'}_{3}\\
  Y^{'}_{4}\\
   \end{array}\right].$$
 As previous,
 assuming  $$\bar{g}_{\beta}(t)=\left[  \begin{array}{cccc}
  \bar{A} & 0 & 0 & 0 \\
   0 & \bar{B} & 0 & 0 \\
   0 & 0 & \bar{C} & 0\\
  0 & 0 & 0 & \bar{D} \\
   \end{array}\right],$$
then we get $$\bar{g}_{\alpha}(t)=\left[ \begin {array}{cccc} \bar{A}+{a}^{2}\bar{B}+{e}^{2}\bar{C}&a\bar{B}+e\bar{C}a&e\bar{C}&b\bar{A}+a\bar{B}c+e\bar{C}d
\\ \noalign{\medskip}a\bar{B}+e\bar{C}a&\bar{B}+{a}^{2}\bar{C}&a\bar{C}&c\bar{B}+a\bar{C}d\\ \noalign{\medskip}e
\bar{C}&a\bar{C}&\bar{C}&d\bar{C}\\ \noalign{\medskip}b\bar{A}+a\bar{B}c+e\bar{C}d&c\bar{B}+a\bar{C}d&d\bar{C}&{b}^{2}\bar{A}+{c}^{2}\bar{B}+{
d}^{2}\bar{C}+\bar{D}\end {array} \right].$$

 Thus $\bar{g}\in [g]$ if and only if all terms convergence to 0 in the sum
\begin{align*}|
\bar{g}-g|_{g}^{2}=&\left(\frac{A-\bar{A}-a^{2}\bar{B}-e^{2}\bar{C}}{A}\right)^{2}+\left(\frac{B-\bar{B}-a^{2}\bar{C}}{B}\right)^{2}+\left(\frac{C-\bar{C}}{C}\right)^{2}
\\+&\left(\frac{b^{2}\bar{A}+c^{2}\bar{B}+d^{2}\bar{C}+\bar{D}-D}{D}\right)^{2}+\frac{2(a\bar{B}+ae\bar{C})^{2}}{AB}+\frac{2(e\bar{C})^{2}}{AC}\\
+&\frac{2(b\bar{A}+ac\bar{B}+de\bar{C})^{2}}{AD}+\frac{2(a\bar{C})^{2}}{BC}+\frac{2(c\bar{B}+ad\bar{C})^{2}}{BD}+\frac{2(d\bar{C})^{2}}{CD}.\\
\end{align*}
It follows that $\bar{g}\in [g]$ if and only if
$$\frac{\bar{A}}{A}\rightarrow 1,\,\,\frac{\bar{B}}{B}\rightarrow 1,\,\,\frac{\bar{C}}{C}\rightarrow 1,\,\,\frac{\bar{D}}{D}\rightarrow 1.$$

\vspace{1ex}

\noindent {\bf Theorem 7.1.} {\it The class $[g]$  is exactly a 2-parameter family.}

\section {\bf {A7.}  $U3I0.$}
\noindent{\bf (A7i)}  \noindent{\bf 8.1. The class $[g]_{\alpha}$ }
\vspace{1ex}

Recall the results in \cite{IJL06}. The   Ricci flow is
 \begin{align*}
&\frac{dA}{dt} =\frac{B}{C}+\frac{C}{B}+2, &\frac{dB}{dt} =&\frac{C}{A}+\frac{D}{C}-\frac{B^{2}}{AC}, \\
&\frac{dC}{dt} =\frac{B}{A}+\frac{D}{B}-\frac{C^{2}}{AB}, &\frac{dD}{dt} =&-\frac{D^{2}}{BC}.
\end{align*}

The long time behavior of the Ricci flow $g(t)$ is
$$A(t)\sim 4t,\,\,B(t)\sim C(t) \sim \frac{1}{D(t)}\sim t^{1/3}.$$

We also have $D(t)=\lambda_{4}\left(1+3\frac{\lambda_{4}}{\lambda_{2}\lambda_{3}}t\right)^{-1/3}$, $BCD^{2}=\lambda_{2}\lambda_{3}\lambda_{4}^{2}$,
and $AD(B-C)=\lambda_{1}\lambda_{4}(\lambda_{2}-\lambda_{3})$.
\vspace{1ex}

\noindent {\bf Lemma 8.1.} {\it The class $[g]_{\alpha}$ is exactly a 3-parameter family.}
\vspace{1ex}

 \noindent{\bf Proof.} Assume that $$\bar{g}(t)=\bar{A}\omega^{1}\otimes\omega^{1}+\bar{B}\omega^{2}\otimes\omega^{2}+\bar{C}\omega^{3}\otimes\omega^{3}+\bar{D}\omega^{4}\otimes\omega^{4}$$ is another Ricci flow solution with the initial data $(\bar{\lambda}_{1},\bar{\lambda}_{2},\bar{\lambda}_{3},\bar{\lambda}_{4})$.
 Then $\bar{g}\in [g]_{\alpha}$ if and only if all terms converge to $0$ in the sum
$$|g-\bar{g}|^{2}_{g}=\left(\frac{A-\bar{A}}{A}\right)^{2}+\left(\frac{B-\bar{B}}{B}\right)^{2}+\left(\frac{C-\bar{C}}{C}\right)^{2}+\left(\frac{D-\bar{D}}{D}\right)^{2}.$$
By (11) in \cite{IJL06}, it is obvious that $\bar{A}/A $ approaches to 1.
By  the equalities $BCD^{2}=\lambda_{2}\lambda_{3}\lambda_{4}^{2}$ and $AD(B-C)=\lambda_{1}\lambda_{4}(\lambda_{2}-\lambda_{3})$, we have $\lim_{t\rightarrow\infty}\frac{B}{C}=1$. Then $\bar{B}/B\rightarrow 1$, $\bar{C}/C\rightarrow 1$ and $\bar{D}/D\rightarrow 1$ if and only if
$\bar{\lambda}_{2}\bar{\lambda}_{3}\bar{\lambda}_{4}^{2}=\lambda_{2}\lambda_{3}\lambda_{4}^{2}$. In fact, if $\bar{B}/B\rightarrow 1$, $\bar{C}/C\rightarrow 1$ and $\bar{D}/D\rightarrow 1$, then $\frac{\bar{B}\bar{C}\bar{D}^{2}}{BCD^{2}}\rightarrow 1$. This yields $\bar{\lambda}_{2}\bar{\lambda}_{3}\bar{\lambda}_{4}^{2}=\lambda_{2}\lambda_{3}\lambda_{4}^{2}$. If $\bar{\lambda}_{2}\bar{\lambda}_{3}\bar{\lambda}_{4}^{2}=\lambda_{2}\lambda_{3}\lambda_{4}^{2}$  holds, then $\bar{D}/D\rightarrow 1$,
$\lim_{t\rightarrow \infty}\frac{\bar{B}^{2}}{B^2}=\lim_{t\rightarrow \infty}\frac{\bar{B}\bar{C}\cdot(\bar{B}/\bar{C})}{BC(\cdot{B/C})}
=\lim_{t\rightarrow \infty}\frac{\bar{B}\bar{C}}{BC}\cdot \lim_{t\rightarrow \infty}\frac{\bar{B}/\bar{C}}{B/C}=1$ and  $\lim_{t\rightarrow \infty}\frac{\bar{C}^{2}}{C^2}=1$. We can choose $\bar{\lambda}_{1}$, $\bar{\lambda}_{2}$ and $\bar{\lambda}_{4}$ arbitrarily.

\vspace{1ex}

\noindent{\bf 8.2.  $g$ and $\bar{g}$ are diagonal under two different frames $\alpha$ and $\beta$}

\vspace{1ex}
Assume  we can use  $Y_{i}=\Lambda_{i}^{k}X_{k}$ with

 $$\Lambda=\left[  \begin{array}{cccc}
  1 & -a_{3} & a_{1} & a_{6} \\
   0 & 1 & 0 & a_{3} \\
   0 & 0 & 1 & a_{1}\\
  0 & 0 & 0 & 1 \\
   \end{array}\right]$$ and $Y^{'}_{i}=\Lambda_{i}^{'k}X_{k}$ with

  $$\Lambda=\left[  \begin{array}{cccc}
  1 & -a_{3}^{'} & a_{1}^{'} & a_{6}^{'} \\
   0 & 1 & 0 & a_{3}^{'} \\
   0 & 0 & 1 & a_{1}^{'}\\
  0 & 0 & 0 & 1 \\
   \end{array}\right]$$to diagonalize $g_{0}$,  $\bar{g}_{0}$ respectively.

   The transformation for $\mathcal{\beta}$ to $\alpha$ is
   $$\left[  \begin{array}{c}
  Y_{1}\\
   Y_{2} \\
   Y_{3}\\
  Y_{4}\\
   \end{array}\right]=\Lambda\Lambda^{'-1} \left[  \begin{array}{c}
  Y^{'}_{1}\\
   Y^{'}_{2} \\
   Y^{'}_{3}\\
  Y^{'}_{4}\\
   \end{array}\right].$$

    Denoting $$a=-a_{3}+a_{3}^{'},\,\,b=a_{1}-a_{1}^{'},\,\,c=a_{6}-a_{6}^{'}+a_{3}a_{3}^{'}-a_{1}a_{1}^{'}-a_{3}^{'2}+a_{1}^{'2},\,\,
   $$and
   $$A=\left[  \begin{array}{cccc}
  1 & a & b & c \\
   0 & 1 & 0 & -a \\
   0 & 0 & 1 & b\\
  0 & 0 & 0 & 1 \\
   \end{array}\right],$$  we have
   $$\left[  \begin{array}{c}
  Y_{1}\\
   Y_{2} \\
   Y_{3}\\
  Y_{4}\\
   \end{array}\right]=A\left[  \begin{array}{c}
  Y^{'}_{1}\\
   Y^{'}_{2} \\
   Y^{'}_{3}\\
  Y^{'}_{4}\\
   \end{array}\right].$$

As previous,
 assuming  $$\bar{g}_{\beta}(t)=\left[  \begin{array}{cccc}
  \bar{A} & 0 & 0 & 0 \\
   0 & \bar{B} & 0 & 0 \\
   0 & 0 & \bar{C} & 0\\
  0 & 0 & 0 & \bar{D} \\
   \end{array}\right],$$

then we get $$\bar{g}_{\alpha}(t)=\left[ \begin {array}{cccc} \bar{A}+{a}^{2}\bar{B}+{b}^{2}\bar{C}+{c}^{2}\bar{D}&a\bar{B}-ca\bar{D}&b\bar{C}+c\bar{D}
b&c\bar{D}\\ \noalign{\medskip}a\bar{B}-c\bar{D}a&\bar{B}+{a}^{2}\bar{D}&-a\bar{D}b&-a\bar{D}
\\ \noalign{\medskip}b\bar{C}+c\bar{D}b&-a\bar{D}b&\bar{C}+{b}^{2}\bar{D}&b\bar{D}\\ \noalign{\medskip}c\bar{D}&
-a\bar{D}&b\bar{D}&\bar{D}\end {array} \right]
.$$

Thus $\bar{g}\in [g]$ if and only if all terms convergence to 0 in the sum
\begin{align*}|
\bar{g}-g|_{g}^{2}=&\left(\frac{A-\bar{A}-a^{2}\bar{B}-b^{2}\bar{C}-c^{2}\bar{D}}{A}\right)^{2}+\left(\frac{B-\bar{B}-a^{2}\bar{D}}{B}\right)^{2}+\left(\frac{C-\bar{C}-b^{2}\bar{D}}{C}\right)^{2}\\
+&\left(\frac{D-\bar{D}}{D}\right)^{2}+\frac{2(a\bar{B}-ac\bar{D})^{2}}{AB}+\frac{2(b\bar{C}+bc\bar{D})^{2}}{AC}+\frac{2(c\bar{D})^{2}}{AD}\\
+&\frac{2(ab\bar{D})^{2}}{BC}+\frac{2(a\bar{D})^{2}}{BD}+\frac{2(b\bar{D})^{2}}{CD}.\\
\end{align*}
    Then we have $$\frac{\bar{A}}{A}\rightarrow 1,\,\,\frac{\bar{B}}{B}\rightarrow 1,\,\,\frac{\bar{C}}{C}\rightarrow 1,\,\,\frac{\bar{D}}{D}\rightarrow 1.$$
 \vspace{1ex}

\noindent {\bf Theorem 8.1.} {\it The class $[g]$  is exactly a 3-parameter family.}

\hspace{2cm}

\noindent{\bf (A7ii)} \noindent{\bf 8.3. The class $[g]_{\alpha}$ }

\vspace{1ex}

Recall the results in \cite{IJL06}. The Ricci flow equation is

 \begin{align*}
&\frac{dA}{dt} =\frac{B^{2}+2(1+a_{2}^{2})BC+(1-a_{2}^{2})^{2}C^{2}}{BC}, &\frac{dB}{dt} =&\frac{AD-B^{2}+(1-a_{2}^{2})^{2}C^{2}}{AC}, \\
&\frac{dC}{dt} =\frac{AD+B^{2}-(1-a_{2}^{2})^{2}C^{2}}{AB}, &\frac{dD}{dt} =&-\frac{D^{2}}{BC}.
\end{align*}
The solution to the Ricci flow is
 \begin{align*}
&A=\lam_{1}+4t, &B=&(\lam_{2}^{3}+3(1-a_{2}^{2})\lam_{2}\lam_{4}t)^{1/3}, \\
&C=\frac{1}{1-a_{2}^{2}}(\lam_{2}^{3}+3(1-a_{2}^{2})\lam_{2}\lam_{4}t)^{1/3}, &D=&\lam_{2}\lam_{4}(\lam_{2}^{3}+3(1-a_{2}^{2})\lam_{2}\lam_{4}t)^{-1/3}.
\end{align*}

We also have $BD=\lambda_{2}\lambda_{4}$.

\vspace{1ex}

\noindent {\bf Lemma 8.2.} {\it The class $[g]_{\alpha}$ is exactly a 2-parameter family.}
\vspace{1ex}

\noindent{\bf Proof.} Assume $$\bar{g}(t)=\bar{A}\omega^{1}\otimes\omega^{1}+\bar{B}\omega^{2}\otimes\omega^{2}+\bar{C}\omega^{3}\otimes\omega^{3}+\bar{D}\omega^{4}\otimes\omega^{4}$$ is another Ricci flow solution with the initial data $(\bar{\lambda}_{1},\bar{\lambda}_{2},\bar{\lambda}_{3},\bar{\lambda}_{4})$.
 Then $\bar{g}\in [g]_{\alpha}$ if and only if all terms converge to $0$ in the sum
$$|g-\bar{g}|^{2}_{g}=\left(\frac{A-\bar{A}}{A}\right)^{2}+\left(\frac{B-\bar{B}}{B}\right)^{2}+\left(\frac{C-\bar{C}}{C}\right)^{2}+\left(\frac{D-\bar{D}}{D}\right)^{2}.$$
It is easy to see that $\bar{A}/A$ goes to 1.
We also get $\bar{B}/B\rightarrow 1$,  $\bar{C}/C\rightarrow 1$ and  $\bar{D}/D\rightarrow 1$ if and only if  $\bar{\lambda}_{2}\bar{\lambda}_{4}=\lambda_{2}\lambda_{4}$. Once  $\bar{\lambda}_{2}$ is chosen,  $\bar{\lambda}_{4}$ and $\bar{\lambda}_{3}$ will be determined since $\bar{\lambda}_{2}=(1-a_{2}^{2})\bar{\lambda}_{3}$ holds in this calss.
\vspace{1ex}

\noindent{\bf 8.4.  $g$ and $\bar{g}$ are diagonal under two different frames $\alpha$ and $\beta$.}
\vspace{1ex}

Assume that we can use  $Y_{i}=\Lambda_{i}^{k}X_{k}$ with

 $$\Lambda=\left[  \begin{array}{cccc}
  1 & -a_{3} & 0 & a_{6} \\
   0 & 1 & a_{2} & a_{3} \\
   0 & 0 & 1 & 0\\
  0 & 0 & 0 & 1 \\
   \end{array}\right]$$ and $Y^{'}_{i}=\Lambda_{i}^{'k}X_{k}$ with

  $$\Lambda^{'}=\left[  \begin{array}{cccc}
  1 & -a_{3}^{'} & 0 & a_{6}^{'} \\
   0 & 1 & a_{2}^{'}& a_{3}^{'} \\
   0 & 0 & 1 & 0\\
  0 & 0 & 0 & 1 \\
   \end{array}\right]$$to diagonalize $g_{0}$,  $\bar{g}_{0}$ respectively.

   The transformation for $\mathcal{\beta}$ to $\alpha$ is
   $$\left[  \begin{array}{c}
  Y_{1}\\
   Y_{2} \\
   Y_{3}\\
  Y_{4}\\
   \end{array}\right]=\Lambda\Lambda^{'-1} \left[  \begin{array}{c}
  Y^{'}_{1}\\
   Y^{'}_{2} \\
   Y^{'}_{3}\\
  Y^{'}_{4}\\
   \end{array}\right].$$

 Denoting $$a=a_{3}^{'}-a_{3},\,\,b=a_{2}-a_{2}^{'},\,\,c=a_{2}^{'}a_{3}-a_{2}^{'}a_{3}^{'},\,\,d=a_{6}-a_{6}^{'}+a_{3}a_{3}^{'}-a_{3}^{'2}
   $$and
   $$A=\left[  \begin{array}{cccc}
  1 & a & c & d \\
   0 & 1 & b & -a \\
   0 & 0 & 1 & 0\\
  0 & 0 & 0 & 1 \\
   \end{array}\right],$$  we have
   $$\left[  \begin{array}{c}
  Y_{1}\\
   Y_{2} \\
   Y_{3}\\
  Y_{4}\\
   \end{array}\right]=A\left[  \begin{array}{c}
  Y^{'}_{1}\\
   Y^{'}_{2} \\
   Y^{'}_{3}\\
  Y^{'}_{4}\\
   \end{array}\right].$$

As previous,  we get$$\bar{g}_{\alpha}(t)=\left[ \begin {array}{cccc} \bar{A}+{a}^{2}\bar{B}+{c}^{2}\bar{C}+{d}^{2}\bar{D}&a\bar{B}+c\bar{C}b-d\bar{D}a&c
\bar{C}&d\bar{D}\\ \noalign{\medskip}a\bar{B}+c\bar{C}b-d\bar{D}a&\bar{B}+{b}^{2}\bar{C}+{a}^{2}\bar{D}&b\bar{C}&-a\bar{D}
\\ \noalign{\medskip}c\bar{C}&b\bar{C}&\bar{C}&0\\ \noalign{\medskip}d\bar{D}&-a\bar{D}&0&\bar{D}
\end {array} \right].$$

Thus $\bar{g}\in [g]$ if and only if all terms convergence to 0 in the sum
\begin{align*}|
\bar{g}-g|_{g}^{2}=&\left(\frac{A-\bar{A}-a^{2}\bar{B}-c^{2}\bar{C}-d^{2}\bar{D}}{A}\right)^{2}+\left(\frac{B-\bar{B}-b^{2}\bar{C}-a^{2}\bar{D}}{B}\right)^{2}\\
+&\left(\frac{C-\bar{C}}{C}\right)^{2}+\left(\frac{D-\bar{D}}{D}\right)^{2}+\frac{2(a\bar{B}+bc\bar{C}-ad\bar{D})^{2}}{AB}+\frac{2(c\bar{C})^{2}}{AC}\\
+&\frac{2(d\bar{D})^{2}}{AD}+\frac{2(b\bar{C})^{2}}{BC}+\frac{2(a\bar{D})^{2}}{BD}.\\
  \end{align*}
  So we have $b=0$ which implies $a_{2}^{'}=a_{2}$ and  $$\frac{\bar{A}}{A}\rightarrow 1,\,\,\frac{\bar{B}}{B}\rightarrow 1,\,\,\frac{\bar{C}}{C}\rightarrow 1,\,\,\frac{\bar{D}}{D}\rightarrow 1.$$

  Under the condition $a_{2}^{'}=a_{2}$, we get the following theorem.

\vspace{1ex}

 \noindent {\bf Theorem 8.2.} {\it The class $[g]$  is exactly a 2-parameter family.}
~~~~~~~~~~~~~~~~~~~~~~~~~~~~~~~~~~~~~~~~~~~~~~~~~~~~~~~~~~~~~~~~~~~~~~~~~~~~~~~~~~~~~~~~~~~~~~~~~~~~~~~~~~~~~~~~~~~~~~~~~~~~
\section {\bf {A8.}  $U3I2.$}

\noindent{\bf 9.1. The class $[g]_{\alpha}$ }

\vspace{1ex}

Recall the results in \cite{IJL06}. The Ricci flow is
 \begin{align*}
&\frac{dA}{dt} =\frac{C}{B}+\frac{B}{C}-2, &\frac{dB}{dt} =&-\frac{B^{2}}{AC}+\frac{C}{A}+\frac{D}{C}, \\
&\frac{dC}{dt} =-\frac{C^{2}}{AB}+\frac{B}{A}+\frac{D}{B}, &\frac{dD}{dt} =&-\frac{D^{2}}{BC}.
\end{align*}
The long time behavior of the Ricci flow $g(t)$ is $$ A(t)\rightarrow k_{4}/2,\,\,B(t)\rightarrow +\infty,\,\,C(t)\rightarrow +\infty,\,\,D(t)\rightarrow 0^{+},$$
where $ k_{4}=\frac{AD(B+C)}{(BC)^{1/2}D}$. We also have $D(t)=\lambda_{4}\left(1+\frac{3\lambda_{4}}{\lambda_{2}\lambda_{3}}t\right)^{-1/3}$,      $BCD^{2}=\lambda_{2}\lambda_{3}\lambda_{4}^{2}$  and  $AD(B+C)=\lambda_{1}\lambda_{4}(\lambda_{2}+\lambda_{3})$.
\vspace{1ex}

\noindent {\bf Lemma 9.1.} {\it The class $[g]_{\alpha}$ is exactly a 2-parameter family.}
\vspace{1ex}

\noindent {\bf Proof.}
Assume that $$\bar{g}(t)=\bar{A}\omega^{1}\otimes\omega^{1}+\bar{B}\omega^{2}\otimes\omega^{2}+\bar{C}\omega^{3}\otimes\omega^{3}+\bar{D}\omega^{4}\otimes\omega^{4}$$ is another Ricci flow solution with the initial data $(\bar{\lambda}_{1},\bar{\lambda}_{2},\bar{\lambda}_{3},\bar{\lambda}_{4})$.
 Then $\bar{g}\in [g]_{\alpha}$ if and only if all terms converge to $0$ in the sum
$$|g-\bar{g}|^{2}_{g}=\left(\frac{A-\bar{A}}{A}\right)^{2}+\left(\frac{B-\bar{B}}{B}\right)^{2}+\left(\frac{C-\bar{C}}{C}\right)^{2}+\left(\frac{D-\bar{D}}{D}\right)^{2}.$$
It follows from  $BCD^{2}=\lambda_{2}\lambda_{3}\lambda_{4}^{2}$  and  $AD(B+C)=\lambda_{1}\lambda_{4}(\lambda_{2}+\lambda_{3})$ that $\lim_{t\rightarrow \infty}\frac{B}{C}=1$.
 Thus $\bar{A}/A\rightarrow 1$, $\bar{B}/B\rightarrow 1$, $\bar{C}/C\rightarrow 1$ and $\bar{D}/D\rightarrow 1$ if and only if $\bar{\lambda}_{2}\bar{\lambda}_{3}\bar{\lambda}_{4}^{2}=\lambda_{2}\lambda_{3}\lambda_{4}^{2}$ and $\bar{\lambda}_{1}\bar{\lambda}_{4}(\bar{\lambda}_{2}+\bar{\lambda}_{3})=\lambda_{1}\lambda_{4}(\lambda_{2}+\lambda_{3})$.

If  $\bar{\lambda}_{2}$ and $\bar{\lambda}_{4}$  are prescribed, then $\bar{\lambda}_{1}$ and $\bar{\lambda}_{3}$ are determined.

\vspace{1ex}

\noindent{\bf 9.2.  $g$ and $\bar{g}$ are diagonal under two different frames $\alpha$ and $\beta$}
\vspace{1ex}

Assume that we can use  $Y_{i}=\Lambda_{i}^{k}X_{k}$ with

 $$\Lambda=\left[  \begin{array}{cccc}
  1 & a_{3} & a_{1} & a_{6} \\
   0 & 1 & 0 & a_{3} \\
   0 & 0 & 1 & a_{1}\\
  0 & 0 & 0 & 1 \\
   \end{array}\right]$$ and $Y^{'}_{i}=\Lambda_{i}^{'k}X_{k}$ with

  $$\Lambda^{'}=\left[  \begin{array}{cccc}
  1 & a_{3}^{'} & a_{1}^{'} & a_{6}^{'} \\
   0 & 1 & 0 & a_{3}^{'} \\
   0 & 0 & 1 & a_{1}^{'}\\
  0 & 0 & 0 & 1 \\
   \end{array}\right]$$to diagonalize $g_{0}$,  $\bar{g}_{0}$ respectively.

   The transformation for $\mathcal{\beta}$ to $\alpha$ is
   $$\left[  \begin{array}{c}
  Y_{1}\\
   Y_{2} \\
   Y_{3}\\
  Y_{4}\\
   \end{array}\right]=\Lambda\Lambda^{'-1} \left[  \begin{array}{c}
  Y^{'}_{1}\\
   Y^{'}_{2} \\
   Y^{'}_{3}\\
  Y^{'}_{4}\\
   \end{array}\right].$$
 Denoting $$a=a_{3}-a_{3}^{'},\,\,b=a_{1}-a_{1}^{'},\,\,c=a_{6}-a_{6}^{'}-a_{1}a_{1}^{'}-a_{3}a_{3}^{'}+a_{1}^{'2}+a_{3}^{'2}
   $$
   and
   $$A=\left[  \begin{array}{cccc}
  1 & a & b & c \\
   0 & 1 & 0 & a \\
   0 & 0 & 1 & b\\
  0 & 0 & 0 & 1 \\
   \end{array}\right],$$  we have
   $$\left[  \begin{array}{c}
  Y_{1}\\
   Y_{2} \\
   Y_{3}\\
  Y_{4}\\
   \end{array}\right]=A\left[  \begin{array}{c}
  Y^{'}_{1}\\
   Y^{'}_{2} \\
   Y^{'}_{3}\\
  Y^{'}_{4}\\
   \end{array}\right].$$

  As previous, we get$$\bar{g}_{\alpha}(t)=\left[ \begin {array}{cccc} \bar{A}+{a}^{2}\bar{B}+{b}^{2}\bar{C}+{c}^{2}\bar{D}&a\bar{B}+c\bar{D}a&b\bar{C}+c\bar{D}
b&c\bar{D}\\ \noalign{\medskip}a\bar{B}+c\bar{D}a&\bar{B}+{a}^{2}\bar{D}&a\bar{D}b&a\bar{D}\\ \noalign{\medskip}
b\bar{C}+c\bar{D}b&a\bar{D}b&\bar{C}+{b}^{2}\bar{D}&b\bar{D}\\ \noalign{\medskip}c\bar{D}&a\bar{D}&b\bar{D}&\bar{D}\end {array}
 \right]
.$$
Thus $\bar{g}\in [g]$ if and only if all terms convergence to 0 in the sum
\begin{align*}|
\bar{g}-g|_{g}^{2}=&\left(\frac{A-\bar{A}-a^{2}\bar{B}-b^{2}\bar{C}-c^{2}\bar{D}}{A}\right)^{2}+\left(\frac{B-\bar{B}-a^{2}\bar{D}}{B}\right)^{2}+\left(\frac{C-\bar{C}-b^{2}\bar{D}}{C}\right)^{2}\\
+&\left(\frac{D-\bar{D}}{D}\right)^{2}+\frac{2(a\bar{B}+ac\bar{D})^{2}}{AB}+\frac{2(b\bar{C}+bc\bar{D})^{2}}{AC}+\frac{2(c\bar{D})^{2}}{AD}\\
+&\frac{2(ab\bar{D})^{2}}{BC}+\frac{2(a\bar{D})^{2}}{BD}+\frac{2(b\bar{D})^{2}}{CD}.\\
  \end{align*}

Thus we get $a=0$ and $b=0$. Under the condition $a=0$ and $b=0$, it is easy to see that
$$\frac{\bar{A}}{A}\rightarrow 1,\,\,\frac{\bar{B}}{B}\rightarrow 1,\,\,\frac{\bar{C}}{C}\rightarrow 1,\,\,\frac{\bar{D}}{D}\rightarrow 1.$$

\vspace{1ex}

\noindent {\bf Theorem 9.1.} {\it The class $[g]$  is exactly a 2-parameter family.}

\section {\bf {A9.}  $U3S1.$}
\noindent{\bf A9i} In this case, the metric $g(t)$ is a product on $\widehat{SL}(2,\mathbb{R})\times\mathbb{R}$
$$g(t)=g_{SL}(t)+\lam_{4}du^{2}$$ where $g_{SL}(t)=A(t)\omega^{1}\otimes \omega^{1}+B(t)\omega^{2}\otimes \omega^{2}+C(t)\omega^{3}\otimes \omega^{3}$ is a  Ricci flow solution on $\widehat{SL}(2,\mathbb{R})$. Refer to \cite{DK01} for the quasi-convergence on $\widehat{SL}(2,\mathbb{R})$.

\vspace{1ex}

\noindent{\bf A9ii} Recall the results in \cite{IJL06}. The Ricci flow is

 \begin{align*}
&\frac{dA}{dt} =\frac{(B+C)^{2}-A^{2}}{BC}+\frac{-A^{2}+B^{2}}{BD}a_{3}^{2}, &\frac{dB}{dt} =&\frac{(A+C)^{2}-B^{2}}{AC}+\frac{A^{2}-B^{2}}{AD}a_{3}^{2}, \\
&\frac{dC}{dt} =\frac{(A-B)^{2}-C^{2}}{AB}, &\frac{dD}{dt} =&\frac{(A+B)^{2}}{AB}a_{3}^{2}.
\end{align*}
In this case, $A=B$ and the Ricci flow reduces to
$$\frac{dA}{dt}=\frac{C}{A}+2,\,\,\,\,\frac{dC}{dt}=-\frac{C^{2}}{A^{2}},\,\,\,\,\frac{dD}{dt}=4a_{3}^{2}.$$
The long time behavior of the Ricci flow $g(t)$ is $$ A(t)\rightarrow +\infty,\,\,B(t)\rightarrow +\infty,\,\,C(t)\rightarrow \text{constant}>0,\,\,D(t)\rightarrow +\infty.$$

\vspace{1ex}

\noindent{\bf 10.1. The class $[g]_{\alpha}$ }
\vspace{1ex}

\noindent {\bf Lemma 10.1.} {\it The class $[g]_{\alpha}$ is exactly a 2-parameter family.}

\vspace{1ex}

\noindent {\bf Proof.}
Assume that $$\bar{g}(t)=\bar{A}\omega^{1}\otimes\omega^{1}+\bar{B}\omega^{2}\otimes\omega^{2}+\bar{C}\omega^{3}\otimes\omega^{3}+\bar{D}\omega^{4}\otimes\omega^{4}$$ is another Ricci flow solution with the initial data $(\bar{\lambda}_{1},\bar{\lambda}_{2},\bar{\lambda}_{3},\bar{\lambda}_{4})$.
 Then $\bar{g}\in [g]_{\alpha}$ if and only if all terms converge to $0$ in the sum
$$|g-\bar{g}|^{2}_{g}=\left(\frac{A-\bar{A}}{A}\right)^{2}+\left(\frac{B-\bar{B}}{B}\right)^{2}+\left(\frac{C-\bar{C}}{C}\right)^{2}+\left(\frac{D-\bar{D}}{D}\right)^{2}.$$

In fact, we compute $$\lim_{t\rightarrow\infty}\frac{\bar{A}}{A}=\lim_{t\rightarrow\infty}\frac{\frac{d\bar{A}}{dt}}{\frac{dA}{dt}}
=\lim_{t\rightarrow\infty}\frac{\frac{\bar{C}}{\bar{A}}+2}{\frac{C}{A}+2}=1.$$

Similarly, we have $\lim_{t\rightarrow\infty}\frac{\bar{B}}{B}=1$. From the equation of $\bar{D}(t)$ and $D(t)$, we get $\lim_{t\rightarrow\infty}\frac{\bar{D}}{D}=1$.
Then $\bar{g}\in [g]_{\alpha}$ if and only if $\bar{C}$ and $C$ have the same limit. The equations for $A$ and $C$ yield
$\frac{d}{dt}\left(\frac{A+C}{AC^2}\right)=0$. Hence we obtain that $C$ is a function of $A$.  Then the choice of $\bar{\lambda}_{1}$ and $\bar{\lambda}_{4}$ determines a metric in $[g]_{\alpha}$.
\vspace{1ex}

\noindent{\bf 10.2.  $g$ and $\bar{g}$ are diagonal under two different frames $\alpha$ and $\beta$}
\vspace{1ex}

Assume that we can use  $Y_{i}=\Lambda_{i}^{k}X_{k}$ with

 $$\Lambda=\left[  \begin{array}{cccc}
  1 & 0 & 0 & 0 \\
   0 & 1 & 0 & 0 \\
   0 & 0 & 1 & 0\\
  0 & 0 & a_{3} & 1 \\
   \end{array}\right]$$ and $Y^{'}_{i}=\Lambda_{i}^{'k}X_{k}$ with

  $$\Lambda^{'}=\left[  \begin{array}{cccc}
  1 & 0 & 0 & 0 \\
   0 & 1 & 0 & 0 \\
   0 & 0 & 1 & 0\\
  0 & 0 & a_{3}^{'} & 1 \\
   \end{array}\right]$$to diagonalize $g_{0}$,  $\bar{g}_{0}$ respectively.

   The transformation for $\mathcal{\beta}$ to $\alpha$ is
   $$\left[  \begin{array}{c}
  Y_{1}\\
   Y_{2} \\
   Y_{3}\\
  Y_{4}\\
   \end{array}\right]=\Lambda\Lambda^{'-1} \left[  \begin{array}{c}
  Y^{'}_{1}\\
   Y^{'}_{2} \\
   Y^{'}_{3}\\
  Y^{'}_{4}\\
   \end{array}\right].$$
 Denoting $a=a_{3}-a_{3}^{'}
   $
   and
   $$A=\left[  \begin{array}{cccc}
  1 & 0 & 0 & 0 \\
   0 & 1 & 0 & 0 \\
   0 & 0 & 1 & 0\\
  0 & 0 & a & 1 \\
   \end{array}\right],$$  we have
   $$\left[  \begin{array}{c}
  Y_{1}\\
   Y_{2} \\
   Y_{3}\\
  Y_{4}\\
   \end{array}\right]=A\left[  \begin{array}{c}
  Y^{'}_{1}\\
   Y^{'}_{2} \\
   Y^{'}_{3}\\
  Y^{'}_{4}\\
   \end{array}\right].$$

As previous,
 we get$$\bar{g}_{\alpha}(t)=\left[ \begin {array}{cccc} \bar{A}&0&0&0\\ \noalign{\medskip}0&\bar{B}&0&0
\\ \noalign{\medskip}0&0&\bar{C}&a\bar{C}\\ \noalign{\medskip}0&0&a\bar{C}&{a}^{2}\bar{C}+\bar{D}
\end {array} \right].$$
Thus $\bar{g}\in [g]$ if and only if all terms convergence to 0 in the sum
$$|
\bar{g}-g|_{g}^{2}=\left(\frac{A-\bar{A}}{A}\right)^{2}+\left(\frac{B-\bar{B}}{B}\right)^{2}+\left(\frac{C-\bar{C}}{C}\right)^{2}+\left(\frac{D-\bar{D}-a^{2}\bar{C}}{D}\right)^{2}
+\frac{2(a\bar{C})^{2}}{CD}.
 $$
We get $a=0$, and $$\frac{\bar{A}}{A}\rightarrow 1,\,\,\frac{\bar{B}}{B}\rightarrow 1,\,\,\frac{\bar{C}}{C}\rightarrow 1,\,\,\frac{\bar{D}}{D}\rightarrow 1.$$
\vspace{1ex}

\noindent {\bf Theorem 10.2.} {\it The class $[g]$  is exactly a 2-parameter family.}

\newpage

{\bf Acknowledgement}

The work was supported by National Natural Science Foundation of China (Grant No.11001268) and Chinese University Scientific Fund (100192014QJ002). The author would like to thank Professor Xiaodong Cao for his discussion and interest on the problem. The author also thank the referees for their helpful comments and suggestions.


\begin{thebibliography}{10}

\bibitem {CGS09} Xiaodong Cao, John Guckenheimer, Laurent Saloff-Coste. The backward behavior of the Ricci and cross curvature flows on $SL(2,\mathbb{R})$. {\it Comm. Anal. Geom.}, 17(4): 777-796, 2009.

\bibitem {CHL} Xiaodong Cao, Songbo Hou and Laurent Saloff-Coste. Backward Ricci flow on locally homogeneous
4-manifolds. preprint

\bibitem{CNS08}Xiaodong Cao, Yilong Ni, and Laurent Saloff-Coste. Cross curvature flow on locally
homogeneous three-manifolds (I). Pacific J. Math., 236(2):263-281, 2008.

\bibitem {CSC09} Xiaodong Cao and Laurent Saloff-Coste. Backward Ricci flow on locally homogeneous
3-manifolds. {\it Comm. Anal. Geom.}, 17(2):305-325, 2009.

\bibitem {CSC08} Xiaodong Cao and Laurent Saloff-Coste. The cross curvature flow on locally homogeneous
three-manifolds (II). {\it Asian J. Math.}, 13(4): 421-458, 2009.

\bibitem{GD08}Glickenstein, David. Riemannian groupoids and solitons for three-dimensional homogeneous
Ricci and cross-curvature flows. {\it Int. Math. Res. Not. IMRN}, (12):Art. ID rnn034, 49, 2008.

\bibitem{Ham86} Richard S. Hamilton. Four-manifolds with positive curvature operator. {\it J. Diff. Geom.}, 24:153-179, 1986.
\bibitem{Ham97} Richard S. Hamilton. Four-manifolds with positive isotropic curvature. {\it Comm. Anal. Geom.}, 5:1-92, 1997.
\bibitem{HI93} Richard S. Hamilton and James Isenberg. Quasi-convergence of Ricci flow for a class of metrics. {\it Comm. Anal. Geom.},1(4):543-559, 1993.
\bibitem {IJ92}James Isenberg and Martin Jackson. Ricci flow of locally homogeneous geometries on
closed manifolds. {\it J. Diff. Geom.}, 35(3):723-741, 1992.

 \bibitem {IJL06} J. Isenberg, M. Jackson and P. Lu.  Ricci flow on locally homogeneous closed 4-manifolds. { \it Comm.
Anal. Geom.,} 14(2):345-386, 2006.
\bibitem{DK00} Dan Knopf. Quasi-convergence of Ricci flow. {\it Comm. Anal. Geom.}, 8(2):375-391, 2000.
\bibitem{DK01} Dan Knopf and Kevin McLeod. Quasi-convergence of model geometries under the Ricci flow. {\it Comm. Anal. Geom.}, 9(4):879-919, 2001.
\bibitem{JL07} John Lott. On the long-time behavior of type-III Ricci flow solutions. {\it Math. Ann.}, 339(3):627-666, 2007.
\bibitem{MM92} M. MacCallum. On the classification of the real four-dimensional Lie
algebra, in On Einstein's path, ed. by A. Harvey, Springer, 299-317,1992.
\bibitem{JM76}John Milnor. Curvatures of left invariant metrics on Lie groups. {\it Advances in Math.,} 21(3):293-329, 1976.

\end{thebibliography}
 \end{document}